\documentclass{amsart}

      \usepackage{amssymb}
\usepackage[T1]{fontenc}
\usepackage[latin9]{inputenc}
\usepackage{amsmath}
\usepackage{amssymb}
\usepackage[all]{xy}
      \theoremstyle{plain}
      \newtheorem{theorem}{Theorem}[section]
      \newtheorem{lemma}[theorem]{Lemma}

      \theoremstyle{definition}
      \newtheorem{definition}[theorem]{Definition}
      \newtheorem*{example}{Example}
      
      \theoremstyle{remark}
      \newtheorem{remark}[theorem]{Remark}

      \makeatletter
      \def\@setcopyright{}
      \def\serieslogo@{}
      \makeatother
   
\begin{document}

%



   \author{Mehdi Khorami}




   \title{A universal coefficient theorem for twisted $K$-theory}


   \begin{abstract}
     In this paper, we recall the definition of twisted $K$-theory in various settings. We prove that for a twist $\tau$ corresponding to a three dimensional integral cohomology class of a space $X$, there exist a ``universal coefficient" isomorphism \[K_{*}^{\tau}(X)\cong K_{*}(P_{\tau})\otimes_{K_{*}(\mathbb{C}P^{\infty})} \hat{K}_{*}\] where $P_\tau$ is the total space of the principal $\mathbb{C}P^{\infty}$-bundle induced over $X$ by $\tau$ and $\hat K_*$ is obtained form the action of  $\mathbb{C}P^{\infty}$ on $K$-theory.
   \end{abstract}



   
   \thanks{The material in this paper will appear in the author's Wesleyan University Ph.D thesis written under supervision of Mark Hovey.}


   \date{January 15, 2010}


   \maketitle



   \section{Introduction}

   In this introduction, we attempt to summarize various equivalent definitions of  twisted $K$-theory and state our results. Here we are only concerned with twists by a 3-dimensional cohomology class of a space $X$. The most common twists of $K$-theory arise from the action of line bundles on $K$-theory. Traditionally, twisted $K$-theory is defined directly by choosing an explicit model for the zero space $BU\times \mathbb Z$ of the $K$-theory spectrum and for $BU(1)\simeq \mathbb{C}P^{\infty}$ and by modeling the action at the point set level. 
   
   Atiyah \cite {Atiyah69} showed that the space $\mathcal F$ of Fredholm operators on an infinite dimensional complex Hilbert space $\mathcal H$ is a representing space for $K$-theory. Roughly speaking, the unitary group $U(\mathcal H)$ acts on the space $\mathcal F$ of Fredholm operators by conjugation (one has to take care of some important details to this approach, including the choice of right topology on $\mathcal F$, which appears in the work of Atiyah and Segal \cite{A.S-TK}). The twisted $K$-cohomology associated to a principal $PU(\mathcal H)$-bundle $P \to X$ is then defined to be the set of homotopy classes of the section of the bundle \[P\times _{PU(\mathcal H)} \mathcal F\to X\]
   
 Thus one twists   $K$ theory by $PU(\mathcal H)$-bundles over $X$. Isomorphism classes of such bundles are classified by maps $X\to BPU(\mathcal H)$. $BPU(\mathcal H)$ is a model for $K(\mathbb{Z},3)$, and thus twists come from homotopy classes of maps $X\to K(\mathbb{Z},3)$.
   
   In essence the action $PU(\mathcal H)\times \mathcal F\to \mathcal F$ corresponds to the action of the Picard group $\text{Pic}(X)$ on $K$-theory:
  \[\text{Pic}(X)\times K(X)\to K(X)\]
  where $\text{Pic}(X)$ is the subgroup of the multiplicative group of the ring $K(X)$ formed by (classes of) line bundles on a compact space $X$ with tensor product as multiplication. On the spectra level, this action is a map  \[ \Sigma^\infty _+\mathbb{C}P^{\infty}\wedge K\to K\]
  
  Construction of twisted cohomology
theories can also be done in the context of parametrized homotopy theory
as suggested by May\textendash{}Sigurdsson (see \cite{M.S-PHT}).
     
     Twisted $K$-theory can be defined using the language of generalized Thom spectra. As far as we know, this is due to Ando, Blumberg and Gepner (``Twisted $K$-theory and $TMF$", to appear). In \cite{units}, the generalized notion of classical Thom spectra was used to define general twists for an $A_\infty $ ring specrum $R$. Any such spectrum has a space of units $\text{GL}_1R$ which deloops to give a classifying space $\text{BGL}_1R$. 
         
    To a  space $X$ and a map \[\tau:X\to \text{BGL}_1R\] there is an associated $R$-algebra Thom spectrum $X^\tau$. In the category of commutative $S$-algebras of \cite{EKMM}, the homotopy type of the Thom spectrum $X^\tau$ can be detected as follows: The principal $\text{GL}_1R$-bundle induces a bundle on $X$
\[
\xymatrix{\text{GL}_1R\ar[d]\ar[r]& \text{GL}_1R\ar[d]\\
Q_\tau \ar[d]\ar[r]& \text{pt} \ar[d]\\
X\ar[r]^\tau \quad & \text{BGL}_1R}
\]
In the derived homotopy category, there is a natural equivalence of $R$-modules  \[X^\tau\simeq \Sigma^\infty _+ Q_\tau\wedge_{\Sigma^\infty _+ \text{GL}_1R}R\]
We warn the reader that this summary neglects many important categorical issues which have been addressed in \cite{units}.

 The $\tau$-twisted $R$-homology of $X$ is defined to be \[R_n^\tau(X)=\pi_n(X^\tau)\] and the $\tau $-twisted  cohomology of $X$ is defined to be\[R^{n}(X)^\tau=\pi_0 R-\text{mod}(X^\tau,\Sigma ^nR)\] 
 
     Ando, Blumberg and Gepner use the Atiyah- Bott-Shapiro orientation map \[\text{MSpin}^c \to K \] to define twisted $K$-theory for twists $\tau:X\to K(\mathbb{Z},3)$. Let $P_\tau$ be the pull-back by the $\mathbb{C}P^{\infty}$-bundle in the diagram 
 \[
 \xymatrix{P_\tau \ar[d] \ar[r]&EK(\mathbb{Z},3)\ar[d]\\
 X\ar[r]^{\tau \quad} &K(\mathbb{Z},3)}
 \] They show that the Thom spectra whose homotopy calculates the $\tau$-twisted $K$-theory is 
 
 \[X^{\tau}\simeq  \Sigma^\infty _+ P_\tau \wedge_{\Sigma^\infty _+ \mathbb{C}P^{\infty}}K\]
 
 Atiyah and Segal showed in \cite{ASTKC} how one can use the Atiyah-Hirebruch spectral sequence with modified differentials to compute twisted K-theory. C. Douglas \cite{CD.TK} uses the twisted the Rothenberg-Steenrod spectral sequence to compute these groups for simply connected simple Lie groups. 
 
 We'll show that there exist a spectral sequence with $E_2$-term 
  \[\text{Tor }_{K_{*}(\mathbb{C}P^{\infty})}^{s,t}(K_{*}(P_{\tau}),\hat{K}_{*})\Rightarrow K_{*}^{\tau}(X)\] converging to the twisted K-theory $K_{*}^{\tau}(X)$. Here $\hat{K}_*$ denotes $K_*$ with $K_{*}(\mathbb{C}P^{\infty})$-module structure obtained from the action map ${\Sigma^\infty _+\mathbb{C}P^{\infty}}_+\wedge K\to K$. The existence of such a spectral sequence may readily be seen for those who are familiar with the methods of \cite{EKMM} if we write \[\Sigma^\infty _+ P_\tau \wedge_{\Sigma^\infty _+ \mathbb{C}P^{\infty}}K\simeq (\Sigma^\infty _+ P_\tau \wedge K) \underset{\Sigma^\infty _+ \mathbb{C}P^{\infty}\wedge K}\wedge K\]However, we'll outline a hands on construction of this spectral sequence in section 3 (See remark 3.7). 
  
  It's easy to see that $\hat{K}_*$ is not flat as a $K_{*}(\mathbb{C}P^{\infty})$-module. Nevertheless, using techniques similar to the ones in the Landweber exact functor theorem \cite{PL.HP}, we'll prove that in the above spectral sequence, all the higher Tor groups vanish: for $s>0$, 
\[\text{Tor }_{K_{*}(\mathbb{C}P^{\infty})}^{s,t}(K_{*}(P_{\tau}),\hat{K}_{*})=0\]
and so \[K_{*}^{\tau}(X)\cong K_{*}(P_{\tau})\otimes_{K_{*}(\mathbb{C}P^{\infty})} \hat{K}_{*}\] 
Here we don't have an analogue of the Landweber filtration theorem, but this is not an obstruction to our proof. The basic idea is that the input modules $K_{*}(P_{\tau})$ have the property that not only they are $K_{*}(\mathbb{C}P^{\infty})$-modules, they are comodules over $K_*K$. 

Thus one can compute twisted $K$-theory by computing the usual $K$-theory of the total space of such bundles $K_*(P_\tau)$. For example, one can use the traditional Serre-Atiyah-Hirebruch spectral sequence to compute these groups. 

  The idea that this should be so occurred to us when investigating a Conner-Floyd type isomorphism for twisted spin bordism and twisted $K$-theory, i.e. whether the Atiyah, Bott and Shapiro orientation map $\text{MSpin}^c \to K $ induces an isomorphism \[\text{MSpin}^{c,\tau}_* (X)\otimes _{\text{MSpin}^c_*} K_*\to K^\tau_*(X)\]      
  The case of the trivial twist $\tau=0$ is in fact a theorem of Hopkins and Hovey in \cite{MMSP} where they prove that the induced map   \[\text{MSpin}^{c}_* (X)\otimes _{\text{MSpin}^c_*} K_*\to K_*(X)\]   is an isomorphism. A proof for the twisted case will appear in a sequel. 
  
  The modern interest in twisted $K$-theory originates from its applications
in mathematical physics, especially string theory. The role of twisted
$K$-theory in string theory started in the 1990's, when Witten argued
that certain topological invariants (charges) of D-branes take values
in twisted $K$-theory of spacetime. More recently, Freed-Hopkins-Teleman showed in a  celebrated theorem  \cite {DFVA,
FHTTK} that the twisted eqivariant  $K$-theory of a compact, connected, and simply connected Lie group $G$ is related to the Verlinde algebra of such a group. 

  \section{preview}
  In this section we build some necessary background for our proofs in the next section.
  
  On the homotopy level, the action of $\mathbb{C}P^{\infty}$ on the $K$-theory spectrum can be analyzed as follows:  let $\mathcal L$ denote the canonical 
line bundle over $\mathbb{C}P^{\infty}$ and $\mathcal V_{n}$ denote the standard n-dimensional bundle
over $BU(n)$. The product $\mathcal L\otimes \mathcal V_{n}$ is an n-dimensional bundle on $\mathbb{C}P^{\infty}\times BU$, and thus is classified by a map \[\mathbb{C}P^{\infty}\times BU\stackrel{j}\to BU\]This map describes the action and on the spectrum level is a map \[\Sigma ^{\infty}_+\mathbb{C}P^{\infty}_+ \wedge K\stackrel{j}\to K\]

The inclusion $i:\mathbb{C}P^{\infty}\simeq BU(1)\to BU$ classifies the line bundle $\mathcal L$ (or $\mathcal L-1$, if we think of BU as \textbf{$BU\times\{0\}$}).
So, we get the composite \[\mathbb{C}P^{\infty}\times BU\stackrel{i\times1}\to BU\times BU\stackrel{\varepsilon}\to BU\]
where $\varepsilon$ is the map obtained by the tensor product of vector bundles. The following diagram commutes.

\[
\xymatrix{\mathbb{C}P^{\infty}\times BU\ar[d]_{i\times1}\ar[dr]^{j}\\
BU\times BU\ar[r] & BU}
\]
On the spectrum level, the inclusion $i$ above is a map $\Sigma ^{\infty}_+\mathbb{C}P^{\infty}\stackrel{i}\longrightarrow K $. We obtain a commutative diagram of spectra

\[
\xymatrix{\Sigma ^{\infty}_+\mathbb{C}P^{\infty}\wedge K\ar[d]_{i\wedge1}\ar[dr]^{j}\\
K\wedge K\ar[r]^\varepsilon & K}
\]
 and thus the action on $K$ is obtained by the composite \[
\Sigma ^{\infty}_+\mathbb{C}P^{\infty}\wedge K\stackrel{i\wedge1}\to K\wedge K\stackrel{\varepsilon}\to K\]
On homotopy, this is a map of $K_*$ algebras \[ K_*(\mathbb{C}P^{\infty})\to K_*K\to K_*\] that we now proceed to describe.

It is well known that in cohomology the ring $K^*(\mathbb{C}P^{\infty})$ is a power series on one generator $x=\gamma -1$ where $\gamma$ is the Hopf line bundle on $\mathbb{C}P^{\infty}$:
\[K^*(\mathbb{C}P^{\infty})=K^*(\text{pt})[|x|]\] It follows that there are unique elements $\beta_i \in K_*(\mathbb{C}P^{n})$, $1\le i \le n$, such that $<x^k,\beta _i>=\delta _i^k$, $1\le k \le n$. Moreover $\{\beta _1, \beta_2,\cdots\}$ forms a $K_*$-basis for $K_*(\mathbb{C}P^{\infty})$:
\[K_*(\mathbb{C}P^{\infty})=K_*\{\beta _1, \beta_2,\dots\}\]
where $K_*$ is the ring $ \mathbb{Z}[t,t^{-1}]$. For simplicity we put all the $\beta_i$s in degree zero (and thus we require $<t^{-k}x^k,\beta _i>=\delta _i^k$ instead). The multiplicative structure in $K_*(\mathbb{C}P^{\infty})$ induced by the $H$-space structure of $\mathbb{C}P^{\infty}$ can be described in terms of the multiplicative formal group law as follows (see \cite{RW-CB}). Let $[1](s)=s$ and inductively define $[n](s)=[n-1](s)+s+s[n-1](s)$. Let  \[\beta(r)=\Sigma_{i\ge 0} \beta_ir^i\] In the power series ring $K_*(\mathbb{C}P^{\infty})[|s,t|]$, we have
\begin{itemize} 
\item $\beta(s)\beta(t)=\beta(s+t+st)$
\item $\beta(s)^m=\beta([m](s))$
\end{itemize}
This gives a complete description of the multiplicative structure of $K_*(\mathbb{C}P^{\infty})$.

The description of $K_*K$ however, requires some work. We briefly outline the construction here. The interested reader is referred  to (\cite{adam} part II, 13) for further details. Since \[K_{q}(K)\cong \text{dirlim}_{n}\tilde{K}_{q+2n}(BU)\] 
for all $q\in \mathbb{Z}$, $K_{q}(K)$ is torsion-free and so we can describe $K_{*}K$
as the image of the monomorphism \[K_{*}K\to K_{*}K\otimes \mathbb{Q}\cong \mathbb{Q}[u,v,u^{-1},v^{-1}]\]
We recall the following theorem.
\begin{theorem} The image of $K_* K$ in $K_{*}K\otimes\mathbb{Q}$ consists
of those finite Laurent series $f(u,v)$ such that \[
f(t,kt)\in\mathbb{Z}[t,t^{-1},\frac{1}{k}]\]
 for all $k\in\mathbb{Z}-\{0\}$.
\end{theorem}

Of particular interest are the elements \[p_i (u,v)=\frac 1{i!}v(v-u)(v-2u)\cdots(v-(i-1)u) \] in $K_* K$ for $i\ge 1$. The following lemma explains our interest in these elements.
\begin {lemma} $K_* K$ is generated over $\mathbb{Z}[u,u^{-1},v^{-1}]$ by the polynomials $\{p_1,p_2,\cdots\}$.

\end{lemma}
\begin {proof} Suppose $f(u,v)$ is a finite Laurent series in $\mathbb{Q}[u,v,u^{-1},v^{-1}]$ which satisfies the condition \[
f(t,kt)\in\mathbb{Z}[t,t^{-1},\frac{1}{k}]\]
 for all $k\in\mathbb{Z}-\{0\}$. We prove that $f(u,v)$ can be written as a $\mathbb{Z}[u,u^{-1},v^{-1}]$-linear combination of the polynomials \[p_i (u,v)=\frac 1{i!}v(v-u)(v-2u)\cdots(v-(i-1)u) \hspace {.3 in}i\ge 1 \] A proof for this can be found in \cite {switzer} (lemma 17.33). For completion we sketch the proof here. We use the following algebraic fact: 
 
 \vspace{.1 in}
 \noindent Let $u(t)\in \mathbb{Q}[t]$ be such that for each integer $k$, $u(k)\in \mathbb{Z}$. Then $u$ can be written as a $\mathbb{Z}$-linear combination of the binomial polynomials \[P_i(t)=\frac 1{i!}t(t-1)(t-2)\cdots(t-(i-1)) \]
 for $i\ge0$ ($P_0$ is interpreted as 1).
 
 By separating $f(u,v)$ into its homogenous components,  we can assume that $f$ is homogeneous of degree say $r$. Let \[g(u,v)=u^{-r}f(u,v)\]
 Then $g$ is homogeneous of degree 0 and we must have \[g(u,v)=h(u^{-1}v)\hspace{.2 in} \text{for some}\hspace{.2 in} h(w)\in \mathbb{Q}[w,w^{-1}]\]We'll prove that $h(k)\in \mathbb{Z}$ for all $k\in \mathbb Z$. Let $m$ be the highest power to which any prime occurs in the denominator of coefficients of $h$. Multiplying $h$ by a suitable positive power of $w$, it is sufficient to consider the case in which $h$ is divisible by $w$ and $w^m$. In particular, $h(0)=0$. Thus if $p$ is a prime dividing $k$, then $p$ does not occur in the denominator of $h(k)$. Note that if $k\neq0$, then \[h(k)=g(t,kt)\in \mathbb Z [\frac 1{k}]\] If $p$ is a prime not dividing $k$, then $p$ does not occur in the denominator of $h(k)$ since $h(k)\in \mathbb Z [\frac 1{k}]$. Thus $h(k)\in \mathbb{Z}$, and the result follows from the fact mentioned above.
\end {proof}

Recall that $K_*K$ is an Hopf algebra with coproduct \[\psi:K_{*}K\to K_{*}K\otimes_{K_{*}}K_{*}K
\text{,}\] the augmentation \[\varepsilon:K_{*}K\to K_{*}\text{,}\]  and conjugation \[c:K_{*}K\to K_{*}K\text{,}\] 
along with two right and left unit maps $\eta_{R},\eta_{L}:K_{*}\to K_{*}K $ induced by \[K\simeq S\wedge K\stackrel{\eta\wedge 1}\longrightarrow K\wedge K\] \[ K\simeq K\wedge S\stackrel{1 \wedge \eta}\longrightarrow K\wedge K\] where $\eta:S\to K$ is the unit of $K$. These data make $(K_*,K_*K)$ into a Hopf algebroid. The following is well-known (see\cite{adam} part II 13, for example).

\begin{theorem} The structure maps of $K_* K$ are given  by the following formulae:
\begin{itemize} 
\item $\psi(u)=u\otimes 1$, $\psi(v)=1\otimes v$;
\item $\varepsilon(u)=t=\varepsilon(v)$;
\item $c(u)=v, c(v)=u$.
\end{itemize}
\end{theorem}
\noindent Here $u=\eta_L(t)$ and $v=\eta_R(t)$.

Thus the multiplication $\varepsilon :K_*K\to K_*$ maps $p_1$ to $t$ and and $p_i$ to zero for $i>1$. Furthermore, the map $K_*(\mathbb{C}P^{\infty})\to K_*K$ sends $t$ to $u$ and $t^i\beta_i$ to $p_i$ (see \cite{switzer}).  Therefore, under the composite 
 $ K_*(\mathbb{C}P^{\infty})\to K_*K\to K_*$,
 $\beta_1 $ maps to $1$ and $\beta _i$ map to zero for $i>1$. 

Let us denote by $\hat{K}_{*}$ the $K_*(\mathbb{C}P^{\infty})$-module obtained as the image of the composite  \[ K_*(\mathbb{C}P^{\infty})\stackrel{i_*}\to K_*K\stackrel{\varepsilon}\to K_*\] 
 Thus $\hat{K}_{*}$ as a ring is the same as the ring $K_*=\mathbb{Z}[t,t^{-1}]$ with $K_*(\mathbb{C}P^{\infty})$-module structure induced by the composite above.
From this point of view, the $K_*(\mathbb{C}P^{\infty})$-module structure of $K_*$ is thus obtained by collapsing $\mathbb{C}P^{\infty}$ to a point $\mathbb{C}P^{\infty}\to \text{pt}$.

The following lemma will be important.
\begin {lemma}
The map $i_*:K_*(\mathbb{C}P^{\infty})\to K_*K$ is flat as a map of rings.
\end{lemma}
\begin {proof} We'll prove that this map is a localization, and thus is flat. It's also easy to see that this map is injective. We'll show that $K_* K$ is in fact obtained from $K_*(\mathbb{C}P^{\infty})$ by inverting $v$, i.e. \[K_*(\mathbb{C}P^{\infty})[v^{-1}]\cong K_*K\] Since $v$ is invertible in $K_*K$, $i_*$ extends injectively to a map \[K_*(\mathbb{C}P^{\infty})[v^{-1}]\to K_*K\] Lemma 2.2 shows that this map is also surjective and the result follows. 
\end{proof}
\section{the spectral sequence}
As already briefly explained in the introduction, there is a spectral sequence \[\text{Tor }_{K_{*}(\mathbb{C}P^{\infty})}^{s,t}(K_{*}(P_{\tau}),\hat{K}_{*})\Rightarrow K_{*}^{\tau}(X)\]
converging to the twisted K-theory of $X$ with a twist $\tau$. Here $t$ shows grading and s corresponds to higher Tor groups.  In this section, we show that all the higher Tor groups in this spectral sequence vanish: for $s>0$, 
\[\text{Tor }_{K_{*}(\mathbb{C}P^{\infty})}^{s,t}(K_{*}(P_{\tau}),\hat{K}_{*})=0\]

Note that $\hat {K}_*$ is not flat as a $K_*(\mathbb{C}P^{\infty})$-module. In fact, as we saw in the previous section, $\beta_i $ maps to zero under $K_*(\mathbb{C}P^{\infty})\stackrel{\varepsilon\ i_*}\longrightarrow \hat{K}_*$ for $i\ge 2$, but it's easy to check that multiplication by $\beta_i$ is injective on $K_*(\mathbb{C}P^{\infty})$. Nevertheless, flatness is not a necessary condition for these Tor groups to vanish. The key point is that the input modules $K_{*}(P_{\tau})$ have very special properties. For any such $P_\tau$, $K_{*}(P_{\tau})$ is a $K_*(\mathbb{C}P^{\infty})$-module which is comodule over $K_*K$. Study of comodule structures was proved beneficial in the work of Landweber for complex cobordism(see \cite{PL.HP}). However, here we do not have anything like Landweber filtration theorem. In his notes \cite {HM.SG}, H. Miller explains how in some settings, the problem of Landweber exactness can be addressed in more generality. 

Quite generally, let $(L,W)$ be an Hopf algebroid of which $(MU_*,MU_*MU)$ and $(K_*,K_*K )$ are motivating examples. This, in particular, means that we have ring homomorphisms 
\[ \begin{array}{ll}
\mbox{coproduct} & \Delta :W\to W \otimes  _L W  \\
\mbox{right and left units} & \eta_R,\eta_L:L\to W \\
\mbox{augmentation} & W\stackrel{\epsilon}\to L\\
\mbox{conjugation} & c:W\to W\\

\end{array}\]
 satisfying certain conditions. $\eta_R$ and $\eta_L$ determine two different $L$-module structure on $W$, called the right and left modules structures. A (left) comodule M over the Hopf algebriod $(L,W)$ is a (left) $L$-module $M$ together with a (left) $L$-linear coaction map \[\psi_M:M\to W\otimes _LM\] which is counitary and coassociative. It is well-known that for any flat ring spectrum $E$ (a ring spectrum $E$ is flat if $E\wedge E$ is equivalent to a wedge of suspensions of $E$), $(E_*,E_*E)$ is an Hopf algebriod and that for any $X$, $E_*(X)$ is a left comodule over it. We'll work in the category $\mathcal{C}$ of comodules over $(L,W)$. For what appears in this paper, $L$ and $W$ are both commutative.
 
 Now, let $R$ be a ring equipped with a ring homomorphism \[f:L\to R\] We can thus consider the $L$-module structure induced on $R$ through this homomorphism. In the case of the Hopf algebroid $(MU_*,MU_*MU)$, any such homomorphism corresponds to a formal group law over $R$. Any such map induces a functor from the category $\mathcal{C}$ to the category of $L$-modules by tensor product with $R$ which we'll denote by $F_R$:
 \[F_R:M\to R\otimes _LM\]
 
 The following definition is due to Hovey and Strickland and has appeared in \cite{comod}.
 \begin{definition} A map of rings $f:L\to R$ is called ``Landweber exact for the Hopf algebroid $(L,W)$" if  the functor  \[F_R:M\to R\otimes _LM\] is exact from the category $\mathcal{C}$ of comodules over $(L,W)$ to the category of $L$-modules.
 \end{definition}
The map $f:L\to R$  induces a ring map 
\[f\otimes \eta _R:L\to R\otimes_L W\]
given by $a\mapsto 1\otimes \eta_R (a)$ and coincides with the composite
\[L\stackrel{\eta _R}\longrightarrow W\cong L\otimes_LW\to R\otimes _LW\] The following important observation was first made by G. Laures and later was discovered independently by Hovey and Strickland in \cite{comod}.  It says that the condition of Landweber exactness can be expressed in terms of the map $f\otimes \eta _R$.
\begin{lemma} Let $(L,W)$ be an Hopf algebriod such that $W$ is flat as an $L$-module under $\eta_L$, so it is also flat over $L$ under $\eta_R$. Then the map $f:L\to R$ is Landweber exact if and only if the map $f\otimes \eta _R$ is flat.
\end{lemma}
\begin {remark} 
We also have the left induced map $\eta _L\otimes f$ which is the composite 
\[L\stackrel{\eta _L}\longrightarrow W\cong W\otimes_LL\to W\otimes _LR\]
 taking $a$ to $\eta_L(a)\otimes 1$. There is no intrinsic difference between the map $f\otimes \eta _R$ and the map $\eta _L\otimes f$ in terms of flatness since the inversion isomorphism $c:W\to W$ swaps $\eta_L$ and $\eta_R$, i.e. the following diagram commutes 
\[
\xymatrix{L\ar[r]^{f\otimes \eta_R\qquad}\ar[dr]_{\eta_L\otimes f}& R\otimes_LW \ar[d]^{1\otimes c}\\
&W\otimes _LR}
\]

\end{remark}
\begin{proof}
First assume that $f\otimes \eta_R$ is flat and $u:A\to B$ is a monomorphism in $\mathcal C$. The coaction map $\psi_A$ is a split monomorphism of $L$-modules; the splitting is given by $\epsilon\otimes 1$. Hence $u$ is a retract of $W\otimes_Lu$ as a map of $L$-modules. It follows that $R\otimes_L u$ is a retract of $R\otimes_L W\otimes_Lu$ as a map of $R$-modules. Now since $R\otimes_LW$ is flat over $L$, we conclude that $R\otimes_L u$ is a monomorphism, as required. Note that in this proof we did not need the flatness assumption of $W$.

Conversely, let $A\to B$ be a monomorphism of $L$-modules. Since $W$ is flat, \[W\otimes_LA\to W\otimes_LB\] is a monomorphism of comodules over $(L,W)$. Moreover, $f:L\to R$ is Landweber exact, and it follows that  \[R\otimes_L(W\otimes_LA)\to R\otimes_L(W\otimes_LB)\] is also a monomorphism. Thus $f\otimes \eta_R:L\to R\otimes _LW$ is flat.
\end{proof}
\begin{remark} While the map $\eta_R\otimes f $ is often not amenable to computations (depending on $f$), it provides a conceptually convenient way to treat Landweber exactness. As we'll see, this works well for our purpose.
\end{remark}
As we mentioned above, for any principal $\mathbb{C}P^{\infty}$-bundle $P_\tau \to X$, $K_*(P_{\tau} )$ is a $K_*(\mathbb{C}P^{\infty})$-module which is also a comodule over $K_*K$. The action of $\mathbb{C}P^{\infty}$ on the total space $P_\tau $ induces a map \[K_*(\mathbb{C}P^{\infty}\times P_{\tau})\to K_*(P_\tau) \] Since  $K_*(\mathbb{C}P^{\infty})$ is free over $K_*$, we have an isomorphism \[K_*(\mathbb{C}P^{\infty}\times P_{\tau})\cong K_*( \mathbb{C}P^{\infty})\otimes_{K_*} K_*(P_{\tau} )\]
and using this, we get the desired module structure 
\[K_*( \mathbb{C}P^{\infty})\otimes_{K_*} K_*(P_{\tau} )\to K_*(P_{\tau} )\]

The Hopf algebroid naturally associated to this situation is \[(K_*(\mathbb{C}P^{\infty}),K_*K\otimes_{K_*}K_*( \mathbb{C}P^{\infty}) \]
We briefly point out the various structure maps. The two unit maps are induced by \[K\wedge \mathbb{C}P^{\infty} \simeq K\wedge S\wedge \mathbb{C}P^{\infty}\stackrel{1\wedge\eta\wedge 1}\longrightarrow K\wedge K\wedge \mathbb{C}P^{\infty}\] \[ K\wedge \mathbb{C}P^{\infty} \simeq S\wedge K \wedge \mathbb{C}P^{\infty}\stackrel{1 \wedge \eta \wedge 1}\longrightarrow  K\wedge K\wedge \mathbb{C}P^{\infty}\] 
The right action of $K_*$ on $K_*K$ figures in the isomorphism \[\pi_*(K\wedge K\wedge \mathbb{C}P^{\infty})\cong K_*K\otimes _{K_*}K_*(\mathbb{C}P^{\infty})\]Note that the left unit map gives the $K_*K$-comodule structure of $K_*( \mathbb{C}P^{\infty})$, which is a map \[K_*(\mathbb{C}P^{\infty})\stackrel{\eta_L}\longrightarrow K_*K\otimes_{K_*} K_*(\mathbb{C}P^{\infty})\] The augmentation is induced by the map\[K\wedge K\wedge \mathbb{C}P^{\infty}\stackrel{\varepsilon\wedge 1}\longrightarrow K\wedge \mathbb{C}P^{\infty}\]
Therefore, $K_*(P_\tau)$ is a comodule over this Hopf algebroid where the comodule map is induced by the map \[K\wedge P_\tau \simeq K\wedge S \wedge P_\tau \to K\wedge K \wedge P_\tau \]
to give a map\[K_*(P_\tau)\to K_*K\otimes_{K_*}K_*(P_\tau)\] 
The only modification necessary is a scalar extension to write \[K_*K\otimes_{K_*}K_*(P_\tau)\cong (K_*K\otimes_{K_*} K_*(\mathbb{C}P^{\infty}))\otimes_{K_*(\mathbb{C}P^{\infty})}K_*(P_\tau)\] to get a comodule map \[K_*(P_\tau)\to (K_*K\otimes_{K_*} K_*(\mathbb{C}P^{\infty}))\otimes_{K_*(\mathbb{C}P^{\infty})}K_*(P_\tau)\]

 Our goal is to prove that the the map $K_*(\mathbb{C}P^{\infty})\stackrel{\varepsilon\ i_*}\longrightarrow \hat{K}_*$ is Landweber exact for the Hopf algebroid \[(K_*(\mathbb{C}P^{\infty}),K_*K\otimes _{K_*}K_*(\mathbb{C}P^{\infty}))\]
 Let $\mathcal D$ denote the category of comodules over this Hopf algebroid. In the light of the lemma 3.2 and the remark 3.3, it suffices to prove that the composite \[K_*(\mathbb{C}P^{\infty})\stackrel{\eta_L}\to K_*K\otimes_{K_*} K_*(\mathbb{C}P^{\infty})\to \hat{K}_*\underset{K_*(\mathbb{C}P^{\infty})}{\otimes}(K_*K\otimes _{K_*}K_*(\mathbb{C}P^{\infty}))\cong K_*K\] is flat.
For this we need to compute the left unit \[\eta _L:K_*(\mathbb{C}P^{\infty})\to K_*K\otimes_{K_*} K_*(\mathbb{C}P^{\infty})\]
Here we briefly describe this map. The interested reader is referred to \cite{switzer} for further details. Set  \[p'_i (u,v)=\frac 1{(i+1)!}(v-u)(v-2u)\cdots(v-iu)\in K_{2n}K \] and let $P=1+p'_1+p'_2+\cdots$. On $K_*( \mathbb{C}P^{\infty})$, the left unit is given by \[
\eta_{L}(t^{k}\beta_{k})=\sum_{i+j=k}(P^{j})_{2i}\otimes t^{j}\beta_{j}\] Here $(P^{j})_{2i}$ shows the homogeneous term with degree $2i$ in the expansion of the jth power of $P$. 

We can now describe the above composite. The map $K_*(\mathbb{C}P^{\infty})\stackrel{\varepsilon\ i_*}\longrightarrow \hat{K}_*$ takes $t\beta_{1}$
to $t$ and sends $t^{i}\beta_{i}$ to zero for $i\geq2$. $\eta_R(t)=v$ and $(P)_{2(i-1)}$
is just $p'_{i-1}$. Thus the above composite sends $t^{i}\beta_{i}$
to \[vp'_{i-1}=\frac{v(v-u)\cdots(v-(i-1)u)}{i!}=p_{i}\]
Thus this composite coincides with the map $i_*:K_*(\mathbb{C}P^{\infty})\to K_*K$ described in the previous section which is flat by lemma 2.4.

Thus we have proved the following:
\begin{theorem} The map $K_*(\mathbb{C}P^{\infty})\stackrel{\varepsilon\ i_*}\longrightarrow \hat{K}_*$ is Landweber exact for the Hopf algebriod 
\[(K_*(\mathbb{C}P^{\infty}),K_*K\otimes_{K_*} K_*(\mathbb{C}P^{\infty}))\]
\end{theorem}

It remains to prove that the category $\mathcal D$ has enough projective modules to ensure the existence of projective resolution within the category. In fact we only need to ensure the existence of relative projective modules with respect to $\hat{K}_*$ since we are only concerned with Tor groups involving $\hat{K}_*$. 
\begin{lemma} For any $P_\tau$, $K_*(P_\tau)$ possesses a relative projective resolution with respect to $\hat{K}_*$ within the category $\mathcal D$. 
\end{lemma}
\begin{proof} We start with the space $X_{0}=P_\tau$ and set $E_{0}=X_{0}\times \mathbb{C}P^{\infty}$.
Since we are working with spaces with basepoint, we require all actions to be basepoint preserving. The action of $CP^{\infty}$ on $P_\tau$ can be assumed
to preserve the basepoint by replacing it with $(P_\tau)_{+}$ and
defining the action on + to be trivial. 

let $f_{0}:E_{0}\to X_{0}$ to be the action map
$\mathbb{C}P^{\infty}\times P_{\tau}\to P_\tau$. This will be a $\mathbb{C}P^{\infty}$-map
if we define the action of $\mathbb{C}P^{\infty}$ on $E_{0}$ to be trivial
on the $X_{0}$ factor, so that $(x,g)h=(x,gh)$. But this action does not fix the basepoint. We can correct this
problem by taking $E_{0}$ to be quotient of $X_{0}\times \mathbb{C}P^{\infty}$
with $\{x_{0}\}\times \mathbb{C}P^{\infty}$ collapsed to a point. 
We have an isomorphism:\[
\widetilde K_{*}(E_{0})\cong K_{*}(\mathbb{C}P^{\infty})\otimes_{K_*} \widetilde K_{*}(X_{0})\]
Thus $K_{*}(E_{0})$ is a $K_{*}(\mathbb{C}P^{\infty})$-module
(which is a comodule over $K_*K$), but it may not be free simply because $K_{*}(X_{0})$ may contain
torsion. However,  \begin{eqnarray*}
Tor_{K_{*}(\mathbb{C}P^{\infty})}(K_{*}(E_{0}),\hat{K}_{*}) & \cong & Tor_{K_{*}(\mathbb{C}P^{\infty})}(K_{*}(X_{0})\underset{K_{*}}{\otimes}K_{*}(\mathbb{C}P^{\infty}),\hat{K}_{*})\\
 & \cong & Tor_{K_{*}}(K_{*}(X_{0}),\hat{K}_{*})=0\end{eqnarray*}

The induced map $f_{0*}:K_{*}(E_{0})\longrightarrow K_{*}(X_{0})$
is a surjection since it is a retraction with respect to the inclusion
$X_{0}\hookrightarrow E_{0}$, $x\mapsto(x,e)$. Let $X_{1}=Cf_{0}$,
the mapping cone of $f_{0}$, so that we get a cofibration \[
E_{0}\to X_{0}\to X_{1}\]
For basepoint reasons we should take the redueced mapping cone. The
action of $CP^{\infty}$ extends naturally to an action on the mapping
cone since it is the mapping cone of a $CP^{\infty}$-map. Applying
$K_{*}$ thus gives a long exact sequence\[
\cdots \to K_{n}(E_{0})\to K_{n}(X_{0})\to K_{n}(X_{1})\to K_{n-1}(E_{0})\to \cdots\]
and since $f_{0*}$ is surjective, we can break this up into a short exact sequence \[
0\to(\Sigma^{-1} K)_{*}(X_{1})\to K_{*}(E_{0})\to K_{*}(X_{0})\to 0\]

We now iterate the construction to produce a diagram

\[
\xymatrix{E_{0}\ar[d] & E_{1}\ar[d] & E_{2}\ar[d]\\
X_{0}\ar[r] & X_{1}\ar[r] & \pi_{n}(X)\ar[r] & ...}
\]
 with associated short exact sequences \[
0\to(\Sigma^{-1}K)_{*}(X_{p+1})\to K_{*}(E_{p})\to K_{*}(X_{p})\to 0\]

 These can be spliced together as in the following diagram to produce
a resolution for $K_{*}(P_\tau)$.

\[
\text{\tiny{\ensuremath{\xymatrix  @C=.8pc {\cdots \ar[r] & K_{*}(E_{2})\ar[dr]\ar[rr] &  & K_{*}(E_{1})\ar[dr]\ar[rr] &  & K_{*}(E_{0})\ar[rr] &  & K_{*}(X_{0})\ar[rr] &  & 0\\
 &  & \Sigma^{-1} K_{*}(X_{2})\ar[ur]\ar[dr] &  & \Sigma^{-1}K_{*}(X_{1})\ar[ur]\ar[dr]\\
 & 0\ar[ur] &  & 0\ar[ur] &  & 0}
}}}\]
where each $K_*(E_i)$ is a $K_*( \mathbb{C}P^{\infty})$-module which also a comodule over $K_*K$. In addition \[Tor_{K_{*}(\mathbb{C}P^{\infty})}(K_{*}(E_{i}),\hat{K}_{*})=0\]
\end{proof}

\begin{remark} We can use the resolution in the above proof to construct our spectral sequence by applying  $-\wedge_{\mathbb{C}P^{\infty}}K$. There is an induced
map \[E_{p}\wedge_{\mathbb{C}P^{\infty}}K\longrightarrow X_{p}\wedge_{\mathbb{C}P^{\infty}}K\]
with mapping cone $X_{p+1}\wedge_{\mathbb{C}P^{\infty}}K$. The associated
long exact sequence of reduced K-homology may no longer split. However,
we can assemble all these long exact sequences into a staircase diagram bellow. For abbreviation, let $Y_p=E_{p}\wedge_{\mathbb{C}P^{\infty}}K$ and $Z_p=X_{p}\wedge_{\mathbb{C}P^{\infty}}K$.
\[
\xymatrix{& & \ar[d] &  & \ar[d]&\\
\cdots \ar[r] &\pi_{n}(Y_p)\ar[r] & \pi_{n}(Z_p)\ar[d]\ar[r] & \pi_{n-1}(Y_{p-1})\ar[r] & \pi_{n-1}(Z_{p-1})\ar[d] \ar[r]&\cdots\\
\cdots \ar[r] &\pi_{n}(Y_{p+1})\ar[r] & \pi_{n}(Z_{p+1})\ar[d] \ar[r] & \pi_{n-1}(Y_p)\ar[r] & \pi_{n-1}(Z_p)\ar[d] \ar[r]&\cdots \\
& &  &  & &}
\]
 Thus we have a spectral sequence.

Let us set $E_{p,q}^{1}=\pi_{p+q}(E_{p}\wedge_{\mathbb{C}P^{\infty}}K)$. The
differential $d_{1}$ is the composition of two horizental maps in
the diagram above. We have
\[
\pi_{*}(E_{p}\wedge_{\mathbb{C}P^{\infty}}K)\simeq\pi_{*}((X_{p}\times \mathbb{C}P^{\infty})\wedge_{\mathbb{C}P^{\infty}}K)\simeq K_{*}(X_{p})\simeq K_{*}(E_{p})\underset{K_{*}(\mathbb{C}P^{\infty})}{\otimes}\hat K_{*}\]
Under this isomorphism the differential $d_{1}$ corresponds to $f_{p*}\otimes1$
where $f_{p}$ is the composition \[E_{p}\longrightarrow X_{p}\longrightarrow\Sigma E_{p-1}\]
the second map being part of the mapping cone sequence \[E_{p-1}\longrightarrow X_{p-1}\longrightarrow X_{p}\longrightarrow \Sigma E_{p-1}\]

The $E^{2}$ terms are the homology of this resolution after tensoring
over ${K}_{*}(\mathbb{C}P^{\infty})$ with $\hat K_{*}$ i.e. the homology of the complex\[
\cdots\longrightarrow K_{*}(E_{1})\underset{K_{*}(\mathbb{C}P^{\infty})}{\otimes}\hat K_{*}\longrightarrow K_{*}(E_{0})\underset{K_{*}(\mathbb{C}P^{\infty})}{\otimes}\hat K_{*}\longrightarrow0\]
The point is that even though the resolution in the above proof
may not be necessarily a projective resolution, the fact that
\[ Tor_{K_{*}(\mathbb{C}P^{\infty})}(K_{*}(E_{p}),\hat K_{*})=0\] says 
that each term $K_{*}(E_{p})$ is projective
relative to $\hat K_{*}$.
Therefore \[
E^{2}=Tor_{K_{*}(\mathbb{C}P^{\infty})}(K_{*}(P_\tau,\hat K_{*})\]
as desired.
\end{remark}

Theorem 3.5 combined with lemma 3.6 completes the proof of our main theorem:
\begin{theorem} In the spectral sequence above, we have  
\[\text{Tor }_{K_{*}(\mathbb{C}P^{\infty})}^{s,t}(K_{*}(P_{\tau}),\hat{K}_{*})=0\]
for $s>0$, and so \[K_{*}^{\tau}(X)\cong K_{*}(P_{\tau})\otimes_{K_{*}(\mathbb{C}P^{\infty})} \hat{K}_{*}\]
\end {theorem}
\begin{proof} Take a relative projective resolution $Q$ for $K_*(P_\tau)$ in $\mathcal D$ by lemma 3.6. The Tor groups $\text{Tor }_{K_{*}(\mathbb{C}P^{\infty})}^{s,t}(K_{*}(P_{\tau}),\hat K_*)$ are the homology of the complex $\hat K_*\otimes Q$, which is exact by theorem 3.5. Therefore they vanish for $s>0$. 

\end{proof}
Remembering the description of the action map $K_*(\mathbb{C}P^{\infty})\to \hat{K}_*$ from the previous section, we see that the isomorphism in this theorem establishes the twisted $K$-group $K_{*}^{\tau}(X)$
as a quotient of $K_{*}(P_{\tau})$ by the submodule generated by
$\beta_{1}-1$, $\beta_{2}$, $\beta_{3}, \cdots$.
\[K_{*}^{\tau}(X)\cong \frac{K_{*}(P_{\tau})}{<\beta_{1}-1, \beta_{2}, \beta_{3}, \cdots>}\]

We expect that this theorem makes computation in twisted K-theory slightly easier. The main difficulty is in  the computation of $K_{*}(P_{\tau})$. We suggest, for example, that C. Douglas's computation for simply connected compact simple lie groups can be recovered from this formula.  Here we offer some examples where the theorem makes computations simpler than the previously known spectral sequence calculations.

\begin{example} The $\mathbb Z_2$-graded twisted $K$-homology of the 3-sphere $S^3$ with nonzero twist \[n:S^3\to K(\mathbb{Z},3)\] is know to be $\mathbb Z/n\mathbb Z$. Let $P_n $ denote the induced bundle by  ``n". Using the standard Atiyah-Hirzebruch-Serre spectral sequence, one can show \[K_*(P_n)=K_{*}(\mathbb{C}P^{\infty})/n\beta_1\] and thus the $\mathbb Z_2$-graded ``n"-twisted $K$-homology of $S^3$ recovers 
\[K^{(n)}_*(S^3)\cong (K_*(\mathbb{C}P^{\infty})/n\beta_1)\otimes_{K_{*}(\mathbb{C}P^{\infty})} \hat{K}_{*}\cong K_*/n=\mathbb Z /n\mathbb Z\]
\end{example}

\begin{example} Let $X=K(\mathbb{Z},3)$ and take the twist to be the identity map $id:K(\mathbb{Z},3)\to K(\mathbb{Z},3)$. Then $P_{id}$ is contractible and so $K_*(P_{id})\cong K_*$. In \[K_{*}\otimes_{K_{*}(\mathbb{C}P^{\infty})} \hat{K}_{*}\] we have $1\otimes 1=0$ since in $\hat{K}_{*}$, 1 is the image of $\beta_1$, and $\beta_1$ maps to zero in $K_*$. Thus 
\[K_{*}^{id}(K(\mathbb{Z},3))\cong K_{*}\otimes_{K_{*}(\mathbb{C}P^{\infty})} \hat{K}_{*}=0\]
Quite generally, let $n:K(\mathbb{Z},3)\to K(\mathbb{Z},3)$ be any twist of $K(\mathbb{Z},3)$, where ``n" is a nonzero integer. It's easy to show that $P_n$ has homotopy type of a $K(\mathbb{Z}/n,2)$. Anderson and Hodgkin \cite{AH-KEM} show that $\widetilde{K}_*(K(\mathbb{Z}/n,2))=0$.
Therefore, for any nonzero twist ``n",  \[K^{(n)}_*(K(\mathbb{Z},3)=0\]
On the other hand the space $K(\mathbb{Z},3)$ behaves like a point in the twisted setting. 
\end{example}

It may be possible to obtain similar results for other cohomology theories using similar methods. For example  it'd be interesting to explore the analogue of 3.8 for the equivariant case. We also suggest to investigate a similar result for twisted real KO-theory where $\mathbb{C}P^{\infty}$ is replaced by $ \mathbb{R}P^{\infty}$ and twists correspond to elements of $H^2(X, \mathbb Z_2)$.
\vspace{.2 in}

\noindent \textbf{Acknowledgments.} It's a  pleasure to thank my thesis advisor Mark Hovey for all of his assistance. I'd also like to thank Haynes Miller, Mark Behrens and Mike Hill for many helpful conversations. 







\end{document}